\documentclass[a4paper]{amsart} 

\usepackage{tikz-cd,mathtools,cancel,amsmath,amssymb}

\newcommand{\id}{\mathrm{id}}
\newcommand{\Ka}{\overline{K}_1^{\mathrm{alg}}}
\newcommand{\minusa}{\cancel{\mathrm{a}}}
\newcommand{\Aff}{\mathrm{Aff}}

\newtheorem{theorem}{Theorem}[section]
\newtheorem{proposition}[theorem]{Proposition}

\theoremstyle{definition}
\newtheorem{definition}[theorem]{Definition}
\newtheorem{example}[theorem]{Example}

\numberwithin{equation}{section}

\begin{document}

\title{Abstract classification theorems for amenable $C^*$-algebras}

\author[S. White]{Stuart White}
\address{Stuart White, Mathematical Institute, University of Oxford,
  Oxford, OX2 6GG, United Kingdom}
\email{stuart.white@maths.ox.ac.uk}
\thanks{This work was partially supported by EPSRC EP/R025061/1-2.}

\begin{abstract}
Operator algebras are subalgebras of the bounded operators on a Hilbert space. They divide into two classes: $C^*$-algebras and von Neumann algebras according to whether they are required to be closed in the norm or weak-operator topology, respectively. 

In the 1970s Alain Connes identified the appropriate notion of amenabilty for von Neumann algebras, and used it to obtain a deep internal finite dimensional approximation structure for these algebras. This structure is exactly what is needed for classification, and one of many consequences of Connes' theorem is the uniqueness of amenable II$_1$ factors, and later a complete classification of all simple amenable von Neumann algebras acting on separable Hilbert spaces.

The Elliott classification programme aims for comparable structure and classification results for $C^*$-algebras using operator $K$-theory and traces. The definitive unital classification theorem was obtained in 2015. This is a combination of the Kirchberg--Phillips theorem and the large scale activity in the stably finite case by numerous researchers over the previous 15--20 years.  It classifies unital simple separable amenable $C^*$-algebras satisfying two extra hypotheses: a universal coefficient theorem which computes $KK$-theory in terms of $K$-theory and a regularity hypothesis excluding exotic high dimensional behaviour. Today the regularity hypothesis can be described in terms of tensor products ($\mathcal Z$-stability). These hypotheses are abstract, and there are deep tools for verifying the universal coefficient theorem and $\mathcal Z$-stability in examples.

This article describes the unital classification theorem, its history and context, together with the new abstract approach to this result developed in collaboration with Carri\'on, Gabe, Schafhauser and Tikuisis. This method makes a direct connection to the von Neumann algebraic results, and does not need to obtain any kind of approximation structure inside $C^*$-algebras en route to classification. The companion survey \cite{W:Survey2} focuses on the role of the $\mathcal Z$-stability hypothesis, and the associated work on "regularity."
\end{abstract}

\maketitle

\section{Operator algebras}

This survey aims to describe how the classification theorems for simple amenable $C^*$-algebras parallel Connes' celebrated von Neumann algebra classification results from the 1970s.  The first four sections set the context and are written for a non-expert. Section \ref{S5} focuses on the two key hypotheses required for the classification theorem.  The last four sections give a flavour of some of the ideas involved in these results,  and require increasingly more background.

Operator algebras originate in the mathematical foundations of quantum mechanics.  There are two main classes: $C^*$-algebras and von Neumann algebras.  A $C^*$-algebra is a complex Banach algebra $A$ equipped with an involution $^*$ satisfying the fundamental $C^*$-identity, $\|x^*x\|=\|x\|^2$ for all $x\in A$.  This seemingly innocuous formula binds the algebraic and analytic aspects of the definition together. It is the source of a surprising amount of rigidity: the norm on a $C^*$-algebra is unique and determined algebraically through spectral data.  $C^*$-algebras are based on abstraction of the bounded operators on a Hilbert space $\mathcal H$, and so examples arise from Hilbert space representations. For example, the left regular representation $\lambda:G\to\mathcal B(\ell^2(G))$ of a discrete group $G$ is given by $\lambda_g(f)(h):= f(g^{-1}h)$ for $g,h\in G$ and $f:G\to\mathbb C$ in $\ell^2(G)$.  To this we associate the \emph{reduced group $C^*$-algebra} of $G$, $C^*_\lambda(G)$ as the $C^*$-algebra generated by the unitaries $(\lambda_g)_{g\in G}$. Group $C^*$-algebras provide a framework linking operator algebras with group representation theory and harmonic analysis.

Norm convergence of bounded operators is a restrictive condition, and there are a number of finer topologies on $\mathcal B(\mathcal H)$: the strong operator topology of pointwise convergence, the weak operator topology of pointwise weak convergence, and the weak$^*$-topology from the duality with the trace class operators.  Von Neumann algebras are the $^*$-subalgebras of $\mathcal B(\mathcal H)$ which are closed in these topologies (they all give $^*$-subalgebras the same closure).   The group von Neumann algebra $L(G)$ is the von Neumann algebra generated by the image of the left regular representation, so $L(G)$ is the weak operator closure of $C^*_\lambda(G)$.

The distinction between the norm and strong or weak operator topologies, creates a profound difference between the flavour of $C^*$-algebras and von Neumann algebras.  The former are topological in nature, while the latter are measure theoretic.  This is evident in the abelian case: via the Gelfand transform, every commutative $C^*$-algebra arises as the algebra $C_0(X)$ of continuous functions vanishing at infinity on a locally compact Hausdorff space $X$. While uniform limits of continuous functions are continuous, pointwise limits are only guaranteed to be measurable, and a general abelian von Neumann algebra is $L^\infty(X,\mu)$ for some measure space $(X,\mu)$.  In the setting of discrete abelian groups $G$, these spaces come from the Fourier transform: $C^*(G)\cong C_0(\widehat{G})$ and $L(G)\cong L^\infty(\widehat{G})$, where $\widehat{G}$ is the Fourier dual group.  Considering $G=\mathbb Z$ and $\mathbb Z^2$, we have $C^*_\lambda(\mathbb Z)\not\cong C^*_\lambda(\mathbb Z^2)$ as the Fourier duals $\mathbb T$ and $\mathbb T^2$ are not homeomorphic. The circle and torus are measurably indistinguishable, so $L(\mathbb Z)\cong L(\mathbb Z^2)$. 

Another important family of examples arises from dynamics.  Given a group action $G\curvearrowright X$ we obtain an induced action $\alpha$ on the relevant abelian operator algebra: $C_0(X)$ or $L^\infty(X)$, according to whether $X$ is topological or measurable. Reminiscent of the semidirect product construction from group theory, the reduced crossed product is a $C^*$-algebra $C_0(X)\rtimes_r G$ or von Neumann algebra $L^\infty(X)\rtimes G$ generated by the relevant abelian algebra $C_0(X)$ or $L^\infty(X)$ and unitaries $(u_g)_{g\in G}$ that implement the induced action: $u_gfu_g^*=\alpha_g(f)$ for $g\in G$ and $f$ in the commutative subalgebra. This construction is set up so that the unitaries $u_g$ generate $C^*_\lambda(G)$ or $L(G)$, respectively.  These operator algebras provide tools for examining actions whose quotient spaces are non-Hausdorff, such as the action $\mathbb Z\curvearrowright\mathbb T$ of rotation by an irrational multiple $\theta$ of $2\pi$ leading to the famous irrational rotation $C^*$-algebras $A_\theta=C(\mathbb T)\rtimes_{r,\theta}\mathbb Z$.  

\section{Projections and approximate finite dimensionality}\label{S2}

Classification results for operator algebras go back to the foundational work of Murray and von Neumann on projections. This relativises the classification of closed subspaces of an infinite-dimensional Hilbert space by their dimension to an operator algebra $A$. Projections $p,q\in A$ are equivalent (written $p\sim q$) if there exists $v\in A$ with $v^*v=p$ and $vv^*=q$; $p$ is sub-equivalent to $q$ ($p\precsim q$) if there exists $q_0\leq q$ with $p\sim q_0$. Akin to Dedekind's definitions for sets, a projection is called infinite if it is equivalent to a proper sub-projection of itself, and finite otherwise. Likewise a unital operator algebra is infinite or finite according to the behaviour of its unit. This theory is particularly clean for von Neumann algebras where $p\precsim q$ and $q\precsim p$ imply $p\sim q$, and is at its crispest in the case of \emph{factors}.  A factor is a von Neumann algebra with trivial centre; these are precisely the simple von Neumann algebras.  Factors are the irreducible building blocks of von Neumann algebras and the set of equivalence classes of projections in a factor is always totally ordered. 

Murray and von Neumann used this total ordering to loosely divide factors into types. The type I factors have a discrete order $\{0,1,2,\dots,n\}$ or $\{0,1,2,\dots\}\cup\{\infty\}$. The only examples are bounded operators on a Hilbert space. The interest begins with the type II factors, where the projections form a continuum.  These subdivide further into II$_1$ factors, where the unit is finite, and II$_\infty$ factors where it is infinite.  After rescaling, the Murray--von Neumann equivalence classes of projections in a II$_1$ factor $\mathcal M$ identify with $[0,1]$ and are determined by the \emph{trace}. That is, there is a unique positive linear functional $\tau$ with $\tau(1_{\mathcal M})=1$ satisfying the trace identity $\tau(xy)=\tau(yx)$. Then projections $p,q\in\mathcal M$ satisfy $p\precsim q$ if and only if $\tau(p)\leq\tau(q)$.  Every II$_\infty$ factor is a von Neumann tensor product of a II$_1$ factor and $\mathcal B(\mathcal H)$, so for many purposes the study of II$_\infty$ factors reduces to that of II$_1$ factors.  Finally, there are the type III factors, which are \emph{purely infinite}: all non-zero projections are infinite (and in the separably acting case equivalent).  

Murray and von Neumann constructed an example of a II$_1$ factor $\mathcal R$ as a suitable completion of the algebraic tensor product of infinitely many copies of the algebra $M_2$ of $2\times 2$ matrices, where the trace comes from extending the product of the normalised traces on the matrix algebras. This led them to isolate a key internal approximation property: \emph{hyperfiniteness}. A von Neumann algebra is \emph{hyperfinite} if finite families of operators can be approximated in strong operator topology by finite dimensional subalgebras (in $\mathcal R$ one uses finite tensor products of copies of $M_2$ to perform these approximations).

 Murray and von Neumann's celebrated uniqueness theorem from \cite{MvN.4} shows that all separably acting hyperfinite II$_1$ factors are isomorphic. In particular, there is nothing special about the $2\times 2$ matrices above. Any choice of matrix size would lead to the same hyperfinite II$_1$ factor, denoted $\mathcal R$. 
 
Examples of II$_1$ factors appear naturally from groups and dynamics.  The von Neumann algebra of an infinite discrete group $G$ is a II$_1$ factor if and only if $G$ is an ICC group: one where all non-trivial elements have infinite conjugacy classes.  When $G$ is an inductive limit of finite groups, $L(G)$ will be hyperfinite. If $G$ is additionally ICC (such as the group $S_\infty$ of all finitely supported permutations of the natural numbers) then it will be isomorphic to $\mathcal R$ by Murray and von Neumann's uniqueness theorem. Approximating an irrational rotation by a rational rotation one gets hyperfiniteness for the II$_1$ factors $L^\infty(\mathbb T)\rtimes\mathbb Z$ associated to irrational rotations.

For $C^*$-algebras, \emph{approximate finite dimensionality} (AF) is defined analogously to hyperfiniteness working with the operator norm in place of the strong topology. This changes things significantly: the infinite $C^*$-tensor products $M_{2^\infty}:=\bigotimes_{n=1}^\infty M_2$ and $M_{3^\infty}:=\bigotimes_{n=1}^\infty M_3$ are not isomorphic. Indeed, there is no unital embedding of $M_2$ into $M_{3^\infty}$, as this would give rise (by a perturbation argument) to a unital embedding of $M_2$ into some $M_{3^n}$, which is impossible.  The difference is that norm-close projections are equivalent, whereas strong-operator-close projections need not be.  This line of reasoning led Glimm to classify all such $C^*$-algebra infinite tensor products of matrices (known as uniformly hyperfinite (UHF) algebras) by the supernatural number consisting of the infinite product of primes appearing in the matrix sizes.  Later, Bratelli examined isomorphisms between separable AF $C^*$-algebras in terms of diagrammatic data for inclusions of finite dimensional $C^*$-algebras, and Elliott gave an algebraic classification which turned out to be in terms of $K$-theory.

Operator $K$-theory is a non-commutative generalisation of Atiyah and Hirzebruch's topological $K$-theory. For a unital $C^*$-algebra $A$, $K_0(A)$ is defined as the Grothendieck group of the abelian semigroup of projections in matrices over $A$ up to equivalence (addition is given by a diagonal direct sum). Elliott's invariant for unital AF-algebras is then $(K_0(A),K_0(A)_+,[1_A]_0)$, where $K_0(A)_+$ is the positive cone arising from projections in matrices over $A$, and $[1_A]_0$ is the class of the unit.  The unital case of Elliott's theorem shows that for separable unital AF $C^*$-algebras $A$ and $B$, any isomorphism $\Phi:K_0(A)\to K_0(B)$ with $\Phi(K_0(A)_+)=K_0(B)_+$ and $\Phi([1_A]_0)=[1_B]_0$ comes from an isomorphism $\phi:A\to B$. 

The $K_0$-functor preserves inductive limits, so one can easily compute $K_0(M_{2^\infty})$ as the dyadic rationals $\{m/2^n:m\in\mathbb Z, n=0,1,\dots\}$ with positivity inherited from the positive reals, and $[1_{M_{2^\infty}}]_0=1$, whereas $K_0(M_{3^\infty})$ is the triadic rationals.  

\section{Connes' theorem}

While Murray and von Neumann's uniqueness theorem is a beautiful result, it is not easy to apply directly; obtaining hyperfiniteness explicitly is out of reach in many examples.  Connes resolved this issue in a landmark abstract characterisation from the 1970s.
\begin{theorem}[{Connes, \cite{C:Ann}}]\label{CT} Let $\mathcal M\subset\mathcal B(\mathcal H)$ be a von Neumann algebra. The following are equivalent:
\begin{enumerate}
\item $\mathcal M$ is injective, i.e. there is a linear contraction $\Phi:\mathcal B(\mathcal H)\to\mathcal M$ with $\Phi(x)=x$ for $x\in\mathcal M$.
\item $\mathcal M$ is semidiscrete.
\item $\mathcal M$ is hyperfinite.
\end{enumerate}
\end{theorem}

Injectivity is so named because it is equivalent to injectivity of $\mathcal M$ in the category of von Neumann algebras and completely positive and contractive (cpc) maps. The definition is independent of the representation of $\mathcal M$.  Semidiscreteness is a finite dimensional approximation property, though a priori of a much weaker nature than hyperfiniteness.  It asks for point weak$^*$ approximations of the identity map by cpc maps factoring through finite dimensional $C^*$-algebras --- such maps preserve the adjoint and order structure (at all levels of matrix amplification), but are not required to preserve the product.   Both injectivity and semidiscreteness are much easier to access than hyperfiniteness in examples; indeed the equivalence of amenablility of a discrete group $G$ with both injectivity and semidiscreteness of $L(G)$ significantly predates Connes' theorem.  Likewise all actions of amenable groups on injective von Neumann algebras produce injective crossed products.

The marriage of Connes' structural theorem with Murray and von Neumann's uniqueness theorem gives a definitive classification result --- there is a unique injective separable II$_1$ factor --- with readily verifiable hypotheses. All the von Neumann algebras associated to countably infinite discrete amenable ICC groups are isomorphic. Moreover, Connes was able to  apply Theorem \ref{CT} to his earlier work on type III factors, giving an almost complete classification of separable injective factors.  The puzzle was completed 10 years later when Haagerup established the uniqueness of the hyperfinite III$_1$ factor \cite{H:Acta}.

\begin{theorem}[Connes, Haagerup, Murray-von Neumann]
There is a complete classification of separable injective factors. 
\end{theorem}

The impact of Theorem \ref{CT} goes well beyond classification results. It shows that hyperfiniteness passes to subalgebras of the hyperfinite II$_1$ factor, which is vital to Jones' theory of subfactors, and directly inspired the classification of amenable equivalence relations by Connes, Feldman and Weiss.  

Connes proved Theorem \ref{CT} in the II$_1$ factor case, and deduced the other cases from this.  A central component of his argument is a deep generalisation of characterisations of amenability in terms of invariant means and F\o{}lner sets to traces on operator algebras. As such it is highly non-constructive, and it remains completely out of reach to describe hyperfiniteness of $L(G)$ explicitly in terms of F\o{}lner sets for an amenable group $G$.   Both Connes' theorem and the techniques involved remain completely instrumental today; for example they play a major role in Popa's deformation-rigidity theory.

Various aspects of Connes' arguments, and the later proofs of Theorem \ref{CT} by Haagerup and Popa have heavily influenced developments for $C^*$-algebras.
Here I highlight two ingredients of Connes' proof for later comparison.

\begin{enumerate}
\item Connes approximates the trace on a separable injective II$_1$ factor $\mathcal M$ externally by a sequence of approximately trace preserving, approximately multiplicative cpc maps (in normalised Hilbert--Schmidt norm) to matrices. Such a sequence is conveniently encoded by an embedding of $\mathcal M$ into the ultrapower $\mathcal R^\omega$ of the hyperfinite II$_1$ factor.  

\item  It is clear that the infinite algebraic tensor product of matrices satisfies $\bigotimes_{n=1}^\infty M_2\cong\bigotimes_{n=1}^\infty M_2\otimes\bigotimes_{n=1}^\infty M_2$. This persists in the weak operator closure used to obtain $\mathcal R$, which becomes idempotent for the von Neumann tensor product: $\mathcal R\cong\mathcal R\otimes\mathcal R$.  This is known as \emph{self-absorption}. A major aspect of Connes' proof is to show that $\mathcal R$ is a tensorial unit for all separable injective II$_1$ factors $\mathcal M$, i.e., $\mathcal M\cong\mathcal M\otimes\mathcal R$. This condition had been previously developed by McDuff, who characterised these \emph{McDuff factors} by the existence of approximately central matrix subalgebras (\cite{McDuff}).
\end{enumerate}

\section{Simple nuclear $C^*$-algebras and the Elliott classification programme}

In a nutshell, the Elliott classification programme aims for a $C^*$-analogue of the Connes, Haagerup, and Murray and von Neumann classification theorem.  Ideally, we seek complete classification results, with abstract hypotheses that are widely verifiable in a range of examples. A fundamental question is: \emph{which $C^*$-algebras should be classified, and by what data}? This is illustrated over the next two sections by the crossed products $C(X)\rtimes_r G$ associated to the action of a discrete group $G$ on a compact Hausdorff space $X$.

Semidiscreteness has a direct analogue for $C^*$-algebras --- the \emph{completely positive approximation property} (CPAP). The only change to the definition is to use the point-norm topology in the $C^*$-setting. The CPAP is in turn equivalent to \emph{nuclearity}, which is characterised by the uniqueness of a $C^*$-norm on the algebraic tensor product $A\odot B$ for all $C^*$-algebras $B$. These conditions give the appropriate notion of amenability for $C^*$-algebras. For instance, a discrete group $G$ is amenable if and only if $C^*_\lambda(G)$ is nuclear and nuclearity is preserved by actions of amenable groups. Moreover, as a surprising application of Connes' theorem, a $C^*$-algebra $A$ has the CPAP if and only if its bidual $A^{**}$ (which is naturally a von Neumann algebra) is semidiscrete.  

The most naive attempt to generalise Connes' theorem fails spectacularly: nuclear $C^*$-algebras will rarely be AF (for example, $C[0,1]$ is nuclear, but certainly not $AF$). Nonetheless there are many situations where nuclear $C^*$-algebras have properties akin to those of injective von Neumann algebras, after making some subtle adjustments to allow for topological phenomena. For example, separable nuclear $C^*$-algebras acting on the same separable Hilbert spaces whose unit balls are close in the Hausdorff metric, are spatially isomorphic (\cite{CSSWW}). This is analogous to existing results for injective von Neumann algebras \cite{C:Invent} --- in both cases amenability is the key hypothesis --- but one has to accept somewhat less strong control on the spatial isomorphism in the $C^*$-setting. 

For group actions, $C(X)\rtimes_rG$ is nuclear whenever $G$ is amenable. One can also define amenable actions of non-amenable groups, such as hyperbolic groups acting on their boundary. Then $C(X)\rtimes_rG$ is nuclear precisely when $G\curvearrowright X$ is amenable.

Within the class of nuclear $C^*$-algebras, the classification programme has mainly focused on simple $C^*$-algebras.  From one point of view, these are analogous to factors --- the simple von Neumann algebras --- but, on the other hand (and unlike von Neumann algebras) we cannot decompose a general $C^*$-algebra in terms of simple algebras. Nevertheless, simple $C^*$-algebras have proved a very fertile ground.  For example, the construction of the Cuntz algebras $\mathcal O_n$ led Cuntz to identify the  $C^*$-analogue of type III factors --- \emph{purely infinite simple $C^*$-algebras}.  These have an abundance of projections (simple purely infinite $C^*$-algebras have real rank zero, and hence are the closed linear span of their projections) and any two non-zero projections $p$ and $q$ satisfy $p\precsim q $ and $q\precsim p$. In a von Neumann factor this would force $p\sim q$, but this need not hold in a $C^*$-algebra.  In fact, the equivalence classes of projections can be very complex: any pair $(G_0,G_1)$ of countable abelian groups can appear as the $K$-theory of a separable simple nuclear purely infinite $C^*$-algebra. Two particularly important examples are $\mathcal O_2$ and $\mathcal O_\infty$ which have $K$-theories $(0,\,0)$ and $(\mathbb Z,\,0)$ respectively.  In the setting of amenable group actions $G\curvearrowright X$ paradoxicality can be used to obtain simple purely infinite nuclear $C^*$-algebras.

The Elliott programme  to classify separable simple nuclear $C^*$-algebras began in earnest after twin breakthroughs in the 1990s.  In the infinite setting, Kirchberg's revolutionary work led to the Kirchberg--Phillips Theorem: a complete $K$-theoretic classification of simple separable nuclear and purely infinite $C^*$-algebras --- now called \emph{Kirchberg algebras} --- satisfying Rosenberg and Schochet's universal coefficient theorem (see Section \ref{S5}).  
\begin{theorem}[Unital Kirchberg--Phillips theorem (\cite{K:ICM,P:Doc}]\label{T4.1}
Let $A$ and $B$ be unital Kirchberg algebras satisfying Rosenberg and Schochet's universal coefficient theorem. Then $A\cong B$ if and only if there is an isomorphism $\alpha_*:K_*(A)\stackrel{\cong}\rightarrow K_*(B)$ with $\alpha_0([1_A]_0)=[1_B]_0$.
\end{theorem}

In the finite setting, Elliott's classification of A$\mathbb T$-algebras (inductive limits of $C^*$-algebras $C(\mathbb T,F)$, where $F$ is finite dimensional) of real rank zero by ordered $K$-theory, combined with the Elliott--Evans theorem that irrational rotation $C^*$-algebras are A$\mathbb T$, to spark substantial work on inductive limit $C^*$-algebras. Thomsen soon realised that ordered $K$-theory alone was not enough to classify all simple A$\mathbb T$ algebras. Traces are also needed. 

As with von Neumann algebras, a trace $\tau$ on a simple $C^*$-algebra $A$ ensures finiteness: if $v^*v\leq vv^*$, then $\tau(vv^*-v^*v)=0$, forcing $vv^*=v^*v$.  Moreover, $M_n(A)$ will be finite for all $n\in\mathbb N$, i.e., $A$ is \emph{stably finite}. A deep theorem of Haagerup shows that stably finite simple nuclear (and, more generally, exact) $C^*$-algebras always admit traces.  So one can detect stable finiteness by traces. However, there is not a direct dichotomy between the finite and infinite: R\o{}rdam produced a finite simple nuclear $C^*$-algebra $A$ which is not stably finite (\cite{R:Acta}).  Such an algebra does not have a von Neumann counterpart.

We write $T(A)$ for the set of traces on $A$. When $A$ is unital, $T(A)$ is compact in the weak$^*$-topology, and convex. Moreover $T(A)$ is a Choquet simplex (this ensures that $T(A)$ can be written as an inverse limit of finite simplices, which is very useful in producing $C^*$-algebras with given trace spaces).  We view $T(A)$ as a family of non-commutative measures, and in many examples it can be determined explicitly.  Indeed, traces on simple crossed products $C(X)\rtimes_rG$ correspond to $G$ invariant measures on $X$.  In this way, the irrational rotation $C^*$-algebras have a unique trace coming from the Lebesgue measure on the circle.

It is often convenient to work with the space $\Aff\,T(A)$ of continuous affine functions $T(A)\to\mathbb R$. By Kadison duality, this is canonically dual to $T(A)$. (One reason for working with $A\mapsto \Aff\, T(A)$ is that makes the entire classification invariant covariant). Every trace induces a well-defined order-preserving map $K_0(A)\to \mathbb R$ mapping $[1_A]_0$ to $1$, and so there is a \emph{pairing} between $K$-theory and traces. At the level of affine functions this is given by 
\begin{equation}
\rho_A:K_0(A)\to\Aff\,T(A),\ \rho_A(x)(\tau)=\tau(x),\quad x\in K_0(A),\ \tau\in T(A).
\end{equation}

Putting these ingredients together, the \emph{Elliott invariant} of a unital $C^*$-algebra $A$ is
\begin{equation}
\mathrm{Ell}(A)=(K_0(A),K_0(A)_+,[1_A]_0,K_1(A),\mathrm{Aff}\,T(A),\rho_A),
\end{equation}
and Elliott conjectured that this should classify all non-elementary simple separable unital nuclear $C^*$-algebras, analogously to the classification of injective factors (\cite{E:ICM}).

For the irrational rotation algebras $A_\theta=C(\mathbb T)\rtimes_{r,\theta}\mathbb Z$, one has $K_*(A_\theta)\cong (\mathbb Z^2,\mathbb Z^2)$. The unique trace $\tau$ embeds $K_0(A)\subset \mathbb R$ by $\tau(m,n)=m+\theta n$ for $m,n\in\mathbb Z$ and here the trace determines the positive cone on $K_0(A)$:  $x\geq 0$ in $K_0(A)$ if and only if $x=0$ or $\tau(x)>0$. A computation then shows $A_\theta\cong A_\phi$ if and only if $\phi\in\pm\theta+\mathbb Z$.   

For irrational rotation algebras both $K_0(A)_+$ and the trace pairing carry the same information. In generality, neither of $K_0(A)_+$ and $(T(A),\rho_A)$  can be recovered from the other. However, in all cases where classification has been obtained $K_0(A)_+$ is determined by $T(A)$ as above.  Thus, although the historical evolution of the Elliott invariant includes $K_0(A)_+$ (this dates back to the classification of AF algebras, whereas traces were only added to the invariant later), with hindsight it is perhaps natural to work with the smaller invariant $KT_u(A)$ consisting of $(K_0(A),[1_A]_0,K_1(A),T(A),\rho_A)$.   I prefer this approach for a number of reasons: it makes it clear that the order structure on $K_0$ does not play an explicit role; the range of $KT_u$ on simple separable unital and nuclear $C^*$-algebras is completely understood (it remains a very challenging problem to determine all possible orders on $K_0(A)_+$ when they are not given by the trace pairing); and $KT_u$ interacts very cleanly with the crucial tensorial absorption condition of $\mathcal Z$-stability (see Section \ref{S5.2} below). I choose to use the Elliott invariant in the main classification theorems, reflecting both the important role $\mathrm{Ell}$ has played historically and the amazing vision Elliott showed in his conjecture, and subsequent programme.

\section{The unital classification theorem}\label{S5}

R\o{}rdam's examples (and precursors by Villadsen) showed that there are highly exotic simple nuclear $C^*$-algebras with phenomena  that have no von Neumann algebraic counterpart. Toms later refined these ideas, constructing simple nuclear $C^*$-algebras which can never be classified by reasonably tractable data \cite{T:Ann}.  Thus we need additional, more subtle, hypotheses to divide the classifiable stably finite simple nuclear $C^*$-algebras from the exotic. By the late 1990s two necessary conditions were known:
\begin{enumerate}
\item[(1)] $A$ satisfies the universal coefficient theorem;
 \item[(2)] $A$ is stable under tensoring by the Jiang-Su algebra $\mathcal Z$, i.e. $A\cong A\otimes \mathcal Z$.
\end{enumerate}

\subsection{The universal coefficient theorem}

Kasparov's bivariant $KK$-theory is one of the fundamental tools in the classification of $C^*$-algebras, providing a tool unifying $K$-theory and extension theory ($K_0(A)\cong KK(\mathbb C,A)$ and $\mathrm{Ext}(A)\cong KK(A,C_0(\mathbb R))$).  In fact, Kirchberg and Phillips both showed that equivalence in $KK$-theory (viewed as a very weak kind of homotopy equivalence) gives rise to an isomorphism for Kirchberg algebras.  Thus the Kirchberg--Phillips theorem has the flavour of a homotopy rigidity result.

\begin{theorem}[Classification of unital Kirchberg Algebras by $KK$-theory]\label{T5.1}
Let $A$ and $B$ be unital Kirchberg algebras. Then $A\cong B$ if and only if there is a $KK$-equivalence $\alpha\in KK(A,B)$ with $\alpha\cdot[1_A]_0=[1_B]_0$.
\end{theorem}

The key tool for computing $KK$-theory is Rosenberg and Schochet's universal coefficient theorem (UCT) from \cite{RS}. The Kasparov product gives a map $KK(A,B)\to\mathrm{Hom}(K_*(A),K_*(B))$, and a $C^*$-algebra $A$ satisfies the UCT when this map fits into a short exact sequence
\begin{equation}
0\to \mathrm{Ext}^1_{\mathbb Z}(K_*(A),K_{*+1}(B))\to KK(A,B)\to\mathrm{Hom}(K_*(A),K_*(B))\to 0
\end{equation}
for all separable $C^*$-algebras $B$. In particular, when both $A$ and $B$ satisfy the UCT, an isomorphism $K_*(A)\cong K_*(B)$ lifts to a $KK$-equivalence; this is how one returns from Theorem \ref{T5.1} to the $K$-theoretic classification of Theorem \ref{T4.1}.

Rosenberg and Schochet established their universal coefficient theorem for abelian $C^*$-algebras, and then showed that the class of nuclear $C^*$-algebras satisfying the UCT is closed under various natural operations (in particular, all $C^*$-inductive limits covered by various classification results lie in the UCT class).  It is a major and rather pressing open problem whether all separable nuclear $C^*$-algebras satisfy the UCT, but for the purposes of applying the classification theorem to concrete examples, this is rarely a difficulty.  Pretty much all separable nuclear $C^*$-algebras that can be explicitly described are known to satisfy the UCT, often through Tu's result (building on Higson and Kasparov's work on the Baum--Connes conjecture) that all $C^*$-algebras associated to amenable groupoids satisfy the UCT.  In particular all crossed products $C(X)\rtimes_rG$ coming from amenable actions satisfy the UCT.

Moreover as the examples realising the invariant for Kirchberg algebras all satisfy the
UCT, the UCT is necessary for a $K$-theoretic classification.

\subsection{$\mathcal Z$-stability}\label{S5.2}

The Cuntz algebra $\mathcal O_\infty$ satisfies $K_*(\mathcal O_\infty)\cong K_*(\mathbb C)=(\mathbb Z,0)$.  Since $\mathcal O_\infty$ satisfies the UCT, one can apply the K\"unneth formula to obtain $K_*(A\otimes\mathcal O_\infty)\cong K_*(A)$ for all separable $A$.  So classification predicts that a Kirchberg algebra $A$ should be isomorphic to $A\otimes\mathcal O_\infty$, i.e., $A$ is \emph{$\mathcal O_\infty$-stable}.  This was confirmed in one of Kirchberg's famous absorption theorems:
\begin{theorem}[Kirchberg's absorption theorems (\cite{K:ICM,KP:Crelle})]
Let $A$ be a separable simple nuclear $C^*$-algebra.  Then:

 (1) $A$ is purely infinite if and only if $A\otimes\mathcal O_\infty\cong A$;

(2) $A\otimes\mathcal O_2\cong \mathcal O_2$.
\end{theorem}

Just as Connes goes via $\mathcal R$-stability of injective II$_1$ factors en-route to hyperfiniteness, the published approaches to the classification of Kirchberg algebras all use the $\mathcal O_\infty$-absorption theorem in a crucial way.   With hindsight $\mathcal O_\infty$-absorption is the key hypothesis which enables the classification of Kirchberg algebras. The absorption theorem then enables classification to be accessed via the more elementary condition of pure infiniteness.  This view point is strengthened by Kirchberg's subsequent classification (by ideal related $KK$-theory) of all separable nuclear $\mathcal O_\infty$-stable $C^*$-algebras.

Returning to stably finite $C^*$-algebras, it is natural to ask \emph{what is the right analogue of the hyperfinite II$_1$ factor}?  A naive first answer might be the CAR-algebra, $M_{2^\infty}$, but this is not canonical. Or one could try the universal UHF-algebra  $\mathcal Q=\bigotimes_{n=2}^\infty M_n$, but this is too big: $M_{2^\infty}$ is a tensorial unit for $\mathcal Q$, not the other way round. No UHF-algebra is a tensorial unit for all its fellows. We want a stably finite unital simple $C^*$-algebra with $K$-theory $(\mathbb Z,\ 0)$ generated by the class of the unit and a unique trace, so that (assuming the K\"unneth formula) it will be a tensorial unit at the level of $K$-theory and traces, just as $\mathcal O_\infty$ is for Kirchberg algebras. One such $C^*$-algebra is $\mathbb C$. Are there others?

In the mid-1990s, Elliott constructed infinite dimensional stably finite simple separable nuclear $C^*$-algebras with arbitrary $K$-theory / trace pairings and so implicitly obtained a $C^*$-algebra with the properties above.  A few years later Jiang and Su tackled this question more explicitly from the view point of tensorial absorption, giving another construction of what we now call the Jiang--Su algebra $\mathcal Z$ (though at the time it would not have been obvious that Elliott and Jiang and Su produced the same algebra) proving $\mathcal Z\cong\mathcal Z\otimes\mathcal Z$. Accordingly, it makes sense to consider \emph{$\mathcal Z$-stable} $C^*$-algebras --- those $A$ for which $A\cong A\otimes\mathcal Z$ --- and, moreover, through the later abstract framework of strongly self-absorbing algebras, Winter showed that $\mathcal Z$-stability is, in a precise sense, the minimal tensorial absorption hypothesis akin to the McDuff property of a II$_1$ factor.   

While $\mathcal Z$ is a little tricky to construct, $\mathcal Z$-stability of a separable $C^*$-algebra can be described without direct reference to $\mathcal Z$ and in a comparable fashion to McDuff's characterisation of $\mathcal R$-stable II$_1$ factors in terms of approximately central matrix subalgebras.  Via Matui and Sato's breakthrough \cite{MS:Acta}, this is particularly clean for a stably finite simple separable nuclear $C^*$-algebra which is $\mathcal Z$-stable when it contains tracially large approximately central cones over matrices.  This is vital in Kerr's approach to detecting $\mathcal Z$-stability of crossed products \cite{K:JEMS} (leading to the recent touchstone result is that all free minimal actions of elementary amenable groups on finite dimensional spaces have $\mathcal Z$-stable crossed product \cite{KN:arXiv}).   I  will describe $\mathcal Z$ and $\mathcal Z$-stability further in the companion survey article \cite{W:Survey2}.

Just as the UCT is necessary for classification, so too is $\mathcal Z$-stability.  Not only are the models realising all  $K$-theory / trace pairings $\mathcal Z$-stable, but tensoring by $\mathcal Z$ acts as the identity on $KT_u$, and, when $A$ is exact, the order on $K_0(A\otimes\mathcal Z)$ is always given by traces.  Moreover, combining work of Kirchberg and R\o{}rdam, for a unital separable nuclear $C^*$-algebra without traces, $\mathcal Z$-stability and $\mathcal O_\infty$-stability are equivalent. In particular $\mathcal Z$-stability is a generalisation of the classifiability hypothesis for Kirchberg algebras.

\subsection{The unital classification theorem, dichotomy, and Toms--Winter regularity}

The combined efforts of large numbers of researchers over close to 30 years have culminated in the definitive classification theorems for simple nuclear $C^*$-algebras, providing the topological counterpart to Connes' theorem. The two subtle hypotheses in the previous subsections are sufficient as well as necessary. We state this in the unital case (and call $C^*$-algebras satisfying the hypothesis of the following theorem \emph{classifiable}).  In particular all crossed products $C(X)\rtimes_rG$ arising from free minimal actions of countable elementary amenable groups on compact metrisable spaces of finite covering dimension are classifiable.

\begin{theorem}[The Unital classification theorem]\label{ClassThm}
Let $A$ and $B$ be unital separable simple nuclear $C^*$-algebras which are $\mathcal Z$-stable and satisfy the UCT.  Then $A$ and $B$ are isomorphic if and only if $\mathrm{Ell}(A)\cong \mathrm{Ell}(B)$. \end{theorem}

As noted above the unital classification theorem is accompanied by a "range of the invariant theorem": any $K$-theory / trace pairing can arise.  This can be used to establish properties of all classifiable $C^*$-algebras through models. For example, all stably finite classifiable $C^*$-algebras are approximately subhomogeneous (with at most $2$-dimensional building blocks); inductive limit structure is a consequence of classification. Recently Li used this to show that all classifiable $C^*$-algebras arise from twisted \'etale groupoids (\cite{L:Invent}). 

Although the unital classification theorem covers both stably finite and purely infinite $C^*$-algebras, the two cases are handled separately.  A beautiful dichotomy theorem of Kirchberg shows that any unital simple exact $C^*$-algebra which is a tensor product of two infinite-dimensional $C^*$-algebras (such as a $\mathcal Z$-stable nuclear $C^*$-algebra) is either purely infinite or stably finite. Moreover the presence or absence of traces decides in which camp a classifiable $C^*$-algebra is found.  The unital classification theorem is then the combination of the Kirchberg--Phillips theorem, and the stably finite case of the unital classification theorem. The rest of the article focuses on the stably finite situation.

The stably finite unital classification theorem was originally obtained in 2015 by combining \cite{GLN1,GLN2,EGLN,TWW} (and with $\mathcal Z$-stability being replaced with the at the time stronger hypothesis of finite nuclear dimension, of which a tiny bit more below).  These in turn build on decades of work --- the unital classification theorem is the collective result of the entire $C^*$-research community --- but in this millennium two names stand out: Lin and Winter.  They drove two major strands of activity in parallel: classification through tracial approximations and regularity through dimension and $\mathcal Z$-stability. The cross fertilisations between these directions have been the source of many of the breakthroughs which have peppered the route to Theorem \ref{ClassThm}.

On the classification side, inspired by Popa's local quantisation technique (a $C^*$-version of the ideas in his proof of injectivity implies hyperfiniteness) Lin introduced the notion of tracial approximations in the early 2000s.  These are a kind of internal approximation by subalgebras whose unit is a projection which is uniformly large in trace.  Weakening the approximation in this way allows more algebras to be reached; the class of tracially AF-algebras is larger than the class of AF-algebras.  Lin and collaborators then massively developed these ideas in a huge body of work culminating in \cite{GLN1,GLN2}.

Meanwhile, Winter and his collaborators developed non-commutative versions of covering dimension (decomposition rank, and then nuclear dimension) for $C^*$-algebras through refined versions of the completely positive approximation property.  This is one of the central concepts in the Toms--Winter conjecture, and combining Winter's $\mathcal Z$-stability theorem from \cite{W:Invent} with the recent \cite{CETWW} (which completes a line of work going back to Matui and Sato's breakthrough \cite{MS:DMJ}) a simple separable unital nuclear non-elementary $C^*$-algebra is $\mathcal Z$-stable if and only if it has finite nuclear dimension.  In the setting of the action of a group $G$ on a finite-dimensional compact metrisable space $X$, one can directly estimate the nuclear dimension of $C(X)\rtimes G$ when $G$ is nilpotent, but one cannot expect direct estimates to work much more generally.  In contrast $\mathcal Z$-stability can be obtained much more generally: now when $G$ is elementary amenable. A detailed discussion of regularity is out of scope here; I will discuss this further in \cite{W:Survey2}.

The two strands come together in a number of landmarks, such as Winter's strategy \cite{W:Crelle} for converting a strong form of classification for UHF-stable $C^*$-algebras to the classification of $\mathcal Z$-stable $C^*$-algebras. The later result is used in the monumental work of Gong, Lin and Niu \cite{GLN1,GLN2} to classify all $\mathcal Z$-stable $C^*$-algebras $A$ with the property that for a UHF-algebra $U$, $A\otimes U$ has a certain $1$-dimensional tracial approximation.  As they show, such algebras exhaust the invariant, so the challenge is to access these very general approximations. This is achieved in \cite{EGLN} using \cite{TWW} by combining finite nuclear dimension and the UCT.

The ideas in the 2015 proof of the unital classification theorem, and the regularity programme are described in more detail in Winter's survey \cite{W:ICM}.   In the rest of this article, I will outline some ingredients in a new abstract and short(er) approach to the stably finite unital classification theorem which is joint work with Carri\'o{}n, Gabe, Schafhauser and Tikuisis (\cite{CGSTW}), which makes more direct contact with von Neumann classification results.

\section{Elliott intertwining: Classifying $C^*$-algebras by classifying maps}\label{S6}

The route towards the unital classification theorem, like many classification results before it, is through a classification of maps together with an Elliott intertwining argument.  This overarching technique goes back to Elliott's classification of AF-algebras, and aspects can even be seen in Murray and von Neumann's uniqueness theorem. We start out by revisiting the classification of AF-algebras via a framework which applies much more generally. 

\subsection{Classifying AF algebras}\label{S6.1}

We consider the classification of countable inductive limits of finite-dimensional algebras by ordered $K_0$, dividing this into three steps.

\smallskip\noindent\textbf{Step 1. Classify maps from finite-dimensional $C^*$-algebras}.  A $C^*$-algebra $B$ has \emph{cancelation} if $K_0$ determines Murray--von Neumann equivalence of projections, i.e. $[p]_0=[q]_0\implies p\sim q$.  This ensures that unital $^*$-homomorphisms from finite-dimensional algebras into  $B$ are classified up to unitary equivalence by ordered $K_0$ and the class of the unit. As ever, classification of maps  consists of two components: \emph{existence} (here that any unital ordered $K_0$-morphism  is realised by a unital $^*$-homomorphism) and \emph{uniqueness} (here two $^*$-homomorphisms are unitarily equivalent if and only if they induce the same map on $K_0$).

\smallskip\noindent\textbf{Step 2. Intertwine to classify maps from inductive limits}. Step 1 can be boosted by taking inductive limits to classify maps from a countable inductive limit $A=\overline{\bigcup_{n=1}^\infty A_n}$ of finite dimensional $C^*$-algebras into a $C^*$-algebra $B$ with cancelation. The invariant remains ordered $K_0$ and the class of the unit, but one can only expect uniqueness up to  \emph{approximate unitary equivalence} ($\phi,\psi:A\to B$ are approximately unitary equivalent ($\phi\approx_{\mathrm{au}}\psi$) when for all finite subsets $\mathcal F\subset  A$ and $\epsilon>0$ there exists a unitary $u\in B$ with $\|u\phi(x)u^*-\psi(x)\|<\epsilon$ for $x\in \mathcal F$). 

That the invariant gives uniqueness of maps $A\to B$ up to approximate unitary equivalence is immediate from the uniqueness up to unitary equivalence of maps $A_n\to B$ in step 1.  Existence needs a (one-sided) Elliott intertwining argument. Given a homomorphism $\Phi:(K_0(A),K_0(A)_+,[1_A]_0)\to (K_0(B),K_0(B)_+,[1_B]_0)$, construct compatible $^*$-homomorphisms $\phi_n:A_n\to B$ inductively. Existence in step 1 gives a map $\tilde{\phi}_n$ implementing $\Phi$ on $K_0(A_n)$. Then uniqueness gives a unitary conjugate $\phi_n$ of $\tilde{\phi}_n$ agreeing with the previously defined $\phi_{n-1}$ on $A_{n-1}$.  The resulting map defined on $\bigcup_{n=1}^\infty A_n$ extends by continuity to a $^*$-homomorphism $\phi$ inducing $\Phi$.  

\smallskip\noindent\textbf{Step 3. Symmetrise assumptions to classify separable AF-algebras}. The following abstract form of Elliott's intertwining argument converts a classification of maps up to approximate unitary equivalence to a classification of algebras.

\begin{proposition}[Elliott's two-sided intertwining argument]\label{2Sided}
Suppose that $A$ and $B$ are separable unital $C^*$-algebras and there are $^*$-homomorphisms $\phi:A\to B$ and $\psi:B\to A$ such that $\psi\circ\phi\approx_{\mathrm{au}}\id_A$ and $\phi\circ\psi\approx_{\mathrm{au}}\id_B$.  Then $A$ and $B$ are isomorphic. Moreover $\phi$ and $\psi$ are approximately unitarily equivalent to mutually inverse isomorphisms.
\end{proposition}

In particular, if a functor $F$ classifies unital maps on a class $\mathcal A$ of separable unital $C^*$-algebras  up to approximate unitary equivalence, then $F$ also classifies $\mathcal A$. Indeed, given $A,B\in\mathcal A$ and an isomorphism $\Phi:F(A)\to F(B)$, the existence component of classification gives $^*$-homomorphisms $\phi:A\to B$ and $\psi:B\to A$ with $F(\phi)=\Phi$ and $F(\psi)=\Phi^{-1}$, and the uniqueness component shows that $\phi$ and $\psi$ satisfy the conditions of Proposition \ref{2Sided}.  

This process classifies those $C^*$-algebras which are both inductive limits of finite dimensional $C^*$-algebras and have cancelation by ordered $K_0$ together with the unit.  But the latter hypothesis is readily seen to be automatic for an AF-algebra, so in fact it classifies countable inductive limits of finite dimensional $C^*$-algebras. 

\subsection{Reducing the unital classification theorem to the classification of approximately multiplicative maps}\label{S6.2}

Let us now return to the general setting, and follow the same 3 step strategy.

\smallskip\noindent\textbf{Step 1. Classify approximately multiplicative maps $A\to B$}. When $A=\overline{\bigcup_{n=1}^\infty A_n}$ is an AF-algebra as in Section \ref{S6.1}, a sequence $(\phi_n)$ of $^*$-homomorphisms $A_n\to B$ can be viewed as an approximately multiplicative map on $A$. Indeed, each $\phi_n$ has a cpc extension to $A$ which is approximately multiplicative in that $\|\phi_n(x)\phi_n(y)-\phi_n(xy)\|\to 0$ for $x,y\in A$. Such approximately multiplicative cpc maps $A\to B$ provide a starting point for classification results in the general setting.

A uniqueness theorem is of the form: for all finite subsets $\mathcal F\subset A$ and $\epsilon>0$, there exists a finite subset $\mathcal G\subset A$ and $\delta>0$ such that any two cpc maps $\phi,\psi:A\to B$ which are $(\mathcal G,\delta)$-approximately multiplicative, and approximately agree on the invariant are approximately unitary equivalent on $\mathcal F$ up to $\epsilon$.  Such statements (and their counterparts for existence) can very quickly become a morass of quantifiers, so it is convenient to use sequence algebras or ultraproducts (just as Connes did in his proof that injectivity implies hyperfiniteness).
\begin{definition}
The \emph{sequence algebra} $B_\infty$ of a $C^*$-algebra $B$ is the quotient $\ell^\infty(B)/c_0(B)$.  It is typical to use representative bounded sequences in $B$ to denote elements of $B_\infty$.
\end{definition}

Reindexing --- the art of turning approximate statements into exact ones ---- is a key feature of sequence algebras and ultraproducts.  For example, when $A$ is a separable $C^*$-algebra, and $B$ a unital $C^*$-algebra, $^*$-homomorphisms $\phi,\psi:A\to B_\infty$ are approximately unitary equivalent if and only if they are unitarily equivalent.  So we aim to classify maps $A\to B_\infty$ up to unitary equivalence. Such a result cleanly encodes a classification of approximately multiplicative maps up to approximate unitary equivalence.

\smallskip
\noindent\textbf{Step 2. Intertwine to classify maps $A\to B$.}  Using separability of $A$ and an intertwining argument we boost the classification
of approximately multiplicative maps to a classification of $^{*}$-homomorphisms
$A\to B$. There is a very clean way to do this through \emph{intertwining through reparameterisations}.   Under very mild conditions (namely that inclusions $B\to B_\infty$ induce an injective map at the level of invariants) this uses classification into $B_\infty$ to show that if $\theta:A\to B_\infty$ looks like it factors through $B$ (i.e it has an invariant factoring through $B$), then it is approximately unitary equivalent to a map that really does factor through $B$.

One can also use this approach for finite von Neumann algebras (suitably adjusted to $\|\cdot\|_2$-approximations and ultrapowers).  Here the classification of maps from finite-dimensional $C^*$-algebras into finite von Neumann algebras $\mathcal N$ by traces (essentially Murray and von Neumann's analysis of projections), gives rise to a classification of maps $\mathcal M\to\mathcal N^\omega$ when $\mathcal M$ is separable and hyperfinite, and $\mathcal N^\omega$ is a tracial ultrapower of $\mathcal N$.  Then a one-sided intertwining classifies maps $\mathcal M\to\mathcal N$, and a two-sided intertwining gives Murray and von Neumann's uniqueness of the hyperfinite II$_1$ factor (avoiding a number of explicit perturbation results).

Moreover, via Connes's theorem we have:
\begin{theorem}\label{T6.3}
Maps from a separable nuclear $C^*$-algebra $A$ to a finite von Neumann algebra $\mathcal N$ are classified up to strong operator approximate unitary equivalence by traces.
\end{theorem}
The point is that any map $A\to\mathcal N$ will factor through a tracial von Neumann completion of $A$ which is injective, so hyperfinite.  Then the previous paragraph applies to give Theorem \ref{T6.3}.  We use this well-known result explicitly in our proof of the unital classification theorem, and it is the only point in the argument where we use some form of internal approximation by subalgebras (namely that the finite part of $A^{**}$ is hyperfinite).

 \smallskip
 \noindent\textbf{Step 3. Symmetrise assumptions to classify $C^*$-algebras} In both steps 1 and 2, the assumptions on the domain $A$ and codomain $B$ are likely to be quite different, and there may also be assumptions on the map (such as nuclearity). Now we symmetrise all the assumptions (requiring that the identity maps on all algebras under consideration satisfy any morphism assumptions) and obtain a classification of algebras using the two-sided Elliott intertwining argument (Proposition \ref{2Sided}).  
 
The upshot of this section is that the unital classification theorem can be expected to follow from a classification of maps $A\to B_\infty$. The rest of the article examines this.

\section{The total invariant for classifying approximate multiplicative maps}\label{S7}

Examples from the 1990s show that $K$-theory and traces are not enough to classify $^*$-homomorphisms.  For example, the tensor flip 
\begin{equation}\label{7.1}
\sigma:\mathcal O_3\otimes\mathcal O_3\to\mathcal O_3\otimes\mathcal O_3;\ x\otimes y\mapsto y\otimes x
\end{equation}
on the Cuntz algebra $\mathcal O_3$ is trivial on $K$-theory. However $\sigma\otimes\mathrm{id}_{\mathcal O_3}$ does not act trivially on $K_0$, and so is not approximately inner.  This section discusses the additional ingredients which must be added to the invariant to obtain uniqueness theorems.

The underlying obstruction behind the example in (\ref{7.1}) is found in $K$-theory with coefficients. Introduced by Schochet, the groups $K_*(A;\mathbb Z/n\mathbb Z)$ fit into a natural six-term exact sequence
\begin{equation}
  \label{7.2}
  \begin{tikzcd}
    K_0(A) \ar[r]
      & K_0(A;\mathbb Z/n\mathbb Z) \ar[r,]
      & K_1(A) \ar[d, "\times n"] \\
    K_0(A) \ar[u, "\times n"]
      & K_1(A;\mathbb Z/n\mathbb Z) \ar[l, ]
      & K_1(A)\ar[l,]
  \end{tikzcd}
 \end{equation}
 and so provide a framework for studying torsion in $K$-theory (namely the kernel of the maps of multiplication by $n$) at the $C^*$-algebraic level.   An efficient way to define $K(A;\mathbb Z/n\mathbb Z)$ is as $K_*(A\otimes C_n)$ for any separable nuclear $C^*$-algebra $C_n$ in the UCT class with $K_*(C_n)\cong (\mathbb Z/n\mathbb Z\ ,0)$, such as $C_n=\mathcal O_{n+1}$. The \emph{total $K$-theory} of $A$, $\underline{K}(A)$ is the combination of $K_*(A)$ and $\bigoplus_{n\geq 2}K_*(A;\mathbb Z/n\mathbb Z)$ together with the natural maps in (\ref{7.2}) (and other natural maps connecting the groups with different coefficients). Each of the groups $K_i(A;\mathbb Z/n\mathbb Z)$ is determined by $K_*(A)$ but in an unnatural fashion.   It is on morphisms $\phi:A\to B$ where $\underline{K}(\phi)$ carries more information than $K_*(\phi)$, e.g., $K_*(\sigma_n)=K_*(\mathrm{id}_{\mathcal O_n})$ while $\underline{K}(\sigma_n)\neq \underline{K}(\mathrm{id}_{\mathcal O_n})$. 

While $KK$-theory determines whether UCT-Kirchberg algebras are isomorphic, it is a little too refined for detecting approximate unitary equivalence: maps $\phi,\psi:A\to B$ with $\phi\approx_{\mathrm{au}}\phi$ can differ in $KK(A,B)$.  R\o{}rdam (in the UCT case) and later Dadarlat (in general) identified a quotient $KL(A,B)$  of $KK(A,B)$ which is constant on approximate unitary equivalence classes of morphisms; for Kirchberg algebras, the converse holds and $KL(A,B)$ determines approximate unitary equivalence. One computes $KL(A,B)$ through Dadarlat and Loring's universal multicoefficient theorem (\cite{DL}): $KL(A,B)\cong \mathrm{Hom}(\underline{K}(A),\underline{K}(B))$ whenever $A$ has the UCT. In this way, total $K$-theory classifies morphisms between UCT-Kirchberg algebras.  This works in much more generality (see, for example, Gabe's retreatment of Kirchberg's $\mathcal O_\infty$-stable classification \cite{Gabe}).

\begin{theorem}
Let $A$ be a separable exact unital $C^*$-algebra satisfying the UCT, and let $B$ be a unital simple $\mathcal O_\infty$-stable $C^*$-algebra. Then unital full nuclear $^*$-homorphisms $A\to B_\infty$ are classified up to approximate unitary equivalence by total $K$-theory (with maps preserving the class of the unit in $K_0$).
\end{theorem}

In the stably finite setting, one needs yet further information. Examples due to Nielsen and Thomsen in the setting of $A\mathbb T$-algebras show the importance of a certain algebraic $K_1$-group.  Given a unital $C^*$-algebra $A$, equip $U_\infty(A)=\bigcup_{n=1}^\infty U(M_n(A))$ with the inductive limit topology.  The map $U_\infty(A)\to K_1(A)$ factors through the abelianisation $U_\infty(A)/DU_\infty(A)$,  where $DU_\infty(A)$ is the group generated by commutators in $U_\infty(A)$, but this  functor is not invariant under approximate unitary equivalence of morphisms.  The solution is to form the \emph{Hausdorffised unitary algebraic $K_1$-group}, $\Ka(A)$ of $A$ as $U_\infty(A)/\overline{DU_\infty(A)}$. We write $\minusa_A:\Ka(A)\to K_1(A)$ for the canonical quotient map.

The group $\Ka(A)$ was systematically studied by Thomsen (\cite{Thomsen}) who used the de la Harpe--Skandalis determinant to relate it to $K$-theory and traces through a natural map $\mathrm{Th}_A:\Aff\ T(A)\to \Ka(A)$ which fits into a sequence
\begin{equation}\label{ThomsenSequence}
\Aff\ T(A)\stackrel{\mathrm{Th}_A}\rightarrow\Ka(A)\stackrel{\minusa_A}\rightarrow K_1(A).
\end{equation}
The kernel of $\mathrm{Th}_A$ is precisely the closure of $\rho_A(K_0(A))$ in $\Aff\,T(A)$, so that $\ker\minusa_A\cong \Aff\,T(A)/\overline{\rho_A(K_0(A))}$. This is a divisible group, so there is a non-canonical splitting
\begin{equation}\label{KASplit}
\Ka(A)\cong K_1(A)\oplus \ker\minusa_A.
\end{equation}
Given a $^*$-homomorphism $\phi:A\to B$, there is no reason why $\Ka(\phi)$ should respect the splittings.  In general, there is a \emph{rotation map} $r_\phi:K_1(A)\to \ker\minusa_B$ so that, with respect to decompositions (\ref{KASplit}), $\Ka(\phi):K_1(A)\oplus \ker\minusa_A\to K_1(B)\oplus \ker\minusa_B$ is given by
\begin{equation}\label{7.4}
\Ka(\phi)=\begin{pmatrix}K_1(\phi)&0\\r_\phi&\Ka(\phi)|_{\ker\minusa_A}\end{pmatrix}.
\end{equation}

\begin{example}
Consider the crossed product $A=(\bigotimes_{-\infty}^\infty\mathcal Z)\rtimes \mathbb Z$, where the action is given by a Bernoulli shift on the tensor product.  The Pimsner--Voiculescu $6$-term exact sequence can be used to calculate $K_*(A)\cong (\mathbb Z,\mathbb Z)$, with $K_1(A)$ being generated by the canonical unitary $u$ implementing the action.  Moreover, $A$ has a unique trace (from the unique trace on $\mathcal Z$). So $\ker\minusa_A\cong\Aff\ T(A)/\overline{\rho_A(K_0(A))}\cong \mathbb R/\mathbb Z\cong \mathbb T$, and $\Ka(A)\cong \mathbb Z\oplus \mathbb T$.   The elements $\{\lambda 1_A:\lambda\in \mathbb T\}$ give representatives of $\ker\minusa_A$ (by an easy de la Harpe-Skandalis determinant calculation).

There are just two automorphisms of $K_*(A)$ fixing $[1_A]_0$: the identity and $(m,n)\mapsto (m,-n)$.  These can be paired with any rotation map $\mathbb Z\to \mathbb T$. The combination of $\mathrm{id}_{K_*(A)}$ with rotation map $\lambda$ is implemented by the automorphism of $A$ fixing $\bigotimes_{-\infty}^\infty\mathcal Z$ and sending $u$ to $\lambda u$, while the automorphism of $A$ which reverses the order of the infinite tensor product and sends $u\mapsto \lambda u^*$ acts as the flip on $K_1(A)$ with rotation map $\lambda$.
\end{example}

It turns out that the extra data described above is now enough for uniqueness. It is formalised in the total invariant.

\begin{definition}\label{DefTotalInv}
The \emph{total invariant} $\underline{K}T_u(A)$ of a unital $C^*$-algebra $A$ consists of $\underline{K}(A)$, $\Aff\,T(A)$, and $\Ka(A)$,
together with all the natural maps between these objects.
\end{definition}

While it is necessary to adjoin $\underline{K}(\cdot)$ and $\Ka(\cdot)$ to Elliott's invariant to obtain uniqueness of morphisms, doing so increases the difficulty of proving the corresponding existence result.  We must now determine exactly which maps between these invariants arise from $^*$-homomorphisms.  In addition to the pairing maps $\rho_\bullet$, the maps $\mathrm{Th}_\bullet$ and $\minusa_\bullet$, it turns out that there are natural maps
\begin{equation}\label{NewPairing}
\zeta^{(n)}_A:K_0(A;\mathbb Z/n\mathbb Z)\to \Ka(A),\quad n\geq 2,
\end{equation}
relating total $K$-theory and $\Ka$. Compatibility with the $\zeta^{(n)}_\bullet$ is an extra obstruction for maps $(\underline{K}(A),\Aff\,T(A),\Ka(A))\to (\underline{K}(B),\Aff\,T(B),\Ka(B))$ to come from a $^*$-homomorphism.  We use the last clause of Definition \ref{DefTotalInv} to regard the maps $\zeta^{(n)}_\bullet$ as part of $\underline{K}T_u$ so that by definition $\underline{K}T_u$-morphisms are compatible with $\zeta^{(n)}_\bullet$.  This completes the total invariant --- no more compatibility requirements are needed for an existence theorem.

The maps $\zeta^{(n)}_A$ are a little fiddly to set up in general (see \cite[Section 3]{CGSTW}, which also sets out how they interact with the other natural maps making up $\underline{K}T_u$), but they are readily identified in straightforward examples. For example, under the identifications $K_0(\mathcal Z;\mathbb Z/n\mathbb Z)\cong\mathbb Z/n\mathbb Z$, and $\Ka(\mathcal Z)\cong \mathbb T$, the maps $\zeta^{(n)}_{\mathcal Z}$ are just the inclusions of the $n$th roots of unity into the circle. Moreover, the maps $\zeta^{(n)}_\bullet$ do not play a role when $K_1(A)$ is torsion free; in this case compatibility with the $\zeta^{(n)}_\bullet$ is automatic from the other compatibility requirements.

Everything is now in place to state a general version of the classification of unital approximate morphisms. Note how the hypotheses found in the unital classification theorem split up amongst the domain, co-domain and morphism in Theorem \ref{ClassificationApprox}.

\begin{theorem}[Stably finite classification of approximately multiplicative maps (\cite{CGSTW}) ]\label{ClassificationApprox} Let $A$ be a separable unital nuclear $C^*$-algebra satisfying the UCT, and let $B$ be a unital simple $\mathcal Z$-stable nuclear $C^*$-algebra with $T(B)\neq\emptyset$. Then the total invariant $\underline{K}T_u$ classifies full unital nuclear maps $A\to B_\infty$ up to unitary equivalence.
\end{theorem}

Once Theorem~\ref{ClassificationApprox} has been established, the Elliott intertwining techniques discussed in Section~\ref{S6} can be used to obtain classification results for algebras. Applying Step~2 of Section~\ref{S6} to  Theorem~\ref{ClassificationApprox} classifies unital nuclear maps
$A\to B$, and then symmetrising assumptions following Step~3 of Section~\ref{S6} classifies the algebras in
the unital classification theorem. But the invariant is
$\underline{K}T_{u}$ not $KT_{u}$ or $\mathrm{Ell}$. So the final ingredient in the unital classification theorem is to extend an isomorphism $KT_u(A)\cong KT_u(B)$ to $\underline{K}T_u(A)\cong \underline{K}T_u(B)$. The extension to $\underline{K}$-theory, and $\Ka(\cdot)$ are purely algebraic results appearing in earlier classification work.  A last little detail is required to correct these extensions and ensure compatibility with the $\zeta^{(n)}_{\bullet}$ when $K_1(A)$ has torsion.  Such an extension is highly non-canonical (and typically far from unique).

\section{Quasidiagonality}\label{S8}

It is easier to construct approximately multiplicative maps (existence in Step \ref{S6.1}) as compared with a $^*$-homomorphism (existence in Step \ref{S6.2}). This is exemplified by contrasting \emph{quasidiagonality} with embeddings into the universal UHF-algebra $\mathcal Q$. Voiculescu showed that quasidiagonality of a $C^*$-algebra $A$ can be viewed as an external approximation property: the existence of approximately multiplicative, approximately isometric cpc maps from $A$ into matrix algebras. When $A$ is separable, these can be packaged into to an embedding of $A$ into $\mathcal Q_\infty$ (this characterises quasidiagonality when $A$ is nuclear).  

Many $C^*$-algebras are quasidiagonal; a deep theorem of Voiculescu shows that quasidiagonality is invariant under homotopy, so that all cones $C_0(0,1]\otimes A$ are quasidiagonal --- a result Kirchberg uses in his $\mathcal O_2$-embedding theorem. On the other hand, as $\mathcal Q$ has a faithful trace, no cone over a simple purely infinite $C^*$-algebra embeds in $\mathcal Q$.  As one cannot model an infinite projection in a matrix algebra, quasidiagonal $C^*$-algebras are stably finite.  The Blackadar--Kirchberg problem asks whether this is the only obstruction for nuclear $C^*$-algebras: are all stably finite nuclear $C^*$-algebras quasidiagonal? This question parallels Connes' important observation that injective II$_1$ factors always embed into $\mathcal R^\omega$. Moreover, the  constructions of stably finite simple separable nuclear $C^*$-algebras which exhaust the Elliott invariant are all quasidiagonal. Finding an abstract source of quasidiagonality is necessary for stably finite classification theorems.  

This was achieved for simple stably finite nuclear $C^*$-algebras with the UCT in the quasidiagonality theorem \cite{TWW}.  The idea is to use traces as a kind of measuring device, by showing that all traces on $A$ are quasidiagonal. One definition of quasidiagonality of $\tau\in T(A)$ is the existence of a sequence $(\theta_n)_n$ of approximately multiplicative cpc maps $A\to\mathcal Q$ with $\tau(x)=\lim_n\tau_{\mathcal Q}(\theta_n(x))$ for all $x\in A$. 

\begin{theorem}[{The quasidiagonality theorem \cite{TWW}}]
Let $A$ be a separable nuclear $C^*$-algebra satisfying the UCT.  Then all faithful amenable traces on $A$ are quasidiagonal.   Accordingly, stably finite simple separable nuclear $C^*$-algebras satisfying the UCT are quasidiagonal.
\end{theorem}

With hindsight the quasidiagonality theorem has turned out to be just the right level of difficulty to isolate and simplify fundamental tools in classification. The original proof was inspired by a stable uniqueness across the interval technique from the tracial approximation approach to classification, and the quasidiagonality theorem is then used to construct tracial approximations from abstract conditions in \cite{EGLN}. A major breakthrough was subsequently made by Schafhauser (\cite{S:Crelle}). Reframing the problem in terms of liftings, he gave a conceptual new proof using $\mathrm{Ext}$-groups.  This idea provides the main framework for our approach to the classification of approximate morphisms (as outlined in the next section).

To sketch Schafhauser's plan, we begin with the trace-kernel extension. At this point is preferable to work with ultrapowers rather than sequence algebras, so let $\omega\in\beta\mathbb N\setminus \mathbb N$ be a free ultrafilter.   We form $\mathcal Q_\omega$ as the quotient of $\ell^\infty (\mathcal Q)$ by those sequences $(x_n)$ with $\lim_{n\to \omega}\|x_n\|=0$; it behaves analogously to $\mathcal Q_\infty$. The ultrapower $\mathcal R^\omega$ is the quotient of $\ell^\infty(\mathcal R)$ by those sequences with $\lim_{n\to\omega}\tau(x_n^*x_n)=0$; the point is that this is a von Neumann algebra, whereas a sequence algebra version is not. Since $\mathcal Q$ is weakly dense in $\mathcal R$, Kaplansky's density theorem gives rise to a surjection of $\mathcal Q_\omega$ onto $\mathcal R^\omega$ with kernel $J$. 

  Given a trace $\tau$ on a separable nuclear $C^*$-algebra $A$, one has an embedding $\theta:A\to\mathcal R^\omega$ realising $\tau$ (by Theorem \ref{T6.3}). Via the Choi--Effros lifting theorem in one direction, and Theorem \ref{T6.3} in the other, $\tau$ is quasidiagonal if and only if $\theta$ lifts to $\tilde{\theta}:A\to\mathcal Q_\omega$,
\begin{equation}
\begin{tikzcd}&&&A\ar[d,"\theta"]\ar[dl,dashrightarrow,swap,"\tilde{\theta}"]&\\0\ar[r]&J\ar[r]&\mathcal Q_\omega\ar[r]&\mathcal R^\omega\ar[r]&0.
\end{tikzcd}
\end{equation}
Forming the pullback extension
\begin{equation}
\begin{tikzcd}\eta:0\ar[r]&J\ar[d,equal]\ar[r]&E\ar[r]\ar[d]&A\ar[d]\ar[r]&0\\{\hphantom{\eta:}0}\ar[r]&J\ar[r]&\mathcal Q_\omega\ar[r]&\mathcal R^\omega\ar[r]&0,
\end{tikzcd}
\end{equation}
liftability of $\theta$ is equivalent to the existence of a $^*$-homomorphism splitting $A\to E$ of $\eta$. 

Extension theory provides the ideal tool for tackling problems of this nature, as $\eta$ induces a class in $\mathrm{Ext}(A,J)$.  However, there is a problem, the trace-kernel ideal $J$ appears somewhat unwieldy.  In particular, it is neither stable nor $\sigma$-unital, which is a deterrent to using $\mathrm{Ext}$.  Schafhauser's key observation is that comparison properties of $\mathcal Q$ ensure that $J$ is \emph{separably stable}: for every separable $C^*$-subalgebra $J_0\subset J$, there is a stable separable $J_1$ with $J_0\subset J_1\subset J$.  With a fair bit of care, this is enough stability to use $\mathrm{Ext}$.

A computation using the UCT and the $K$-theory of $\mathcal R^\omega$ readily shows that $[\eta]=0$ in $\mathrm{Ext}(A,J)$. This does not yet mean that $\eta$ splits, but rather that after adjoining a further trivial extension $\eta_1$ say, the sum $\eta\oplus\eta_1$ splits.  So the final step is to ensure that $\eta$ is \emph{absorbing} so that $\eta\cong \eta\oplus\eta_1$. This is achieved using an abstract Weyl--von Neumann / Voiculescu-type theorem of Elliott and Kucerovsky \cite{EK} which heavily exploits Kirchberg's work for infinite $C^*$-algebras. When $A$ is non-unital,  absorption is a consequence of the faithfulness of $\tau$ via injectivity of $\theta$.  There is an important detail when $A$ (and hence $\theta$) are unital. In this case $\eta$ can never be absorbing and we can only ask for absorption of unital extensions.  The trick is to pass to a non-unital $2\times 2$ matrix amplification to replace $\theta$ by a non-unital map. This de-unitisation idea recurs extensively in the classification of approximate morphisms.

\section{Classification of approximately multiplicative maps}\label{S9}

We end with a brief  discussion of some ingredients in the classification of approximately multiplicative maps.  For the rest of the article, let $A$ and $B$ be as in Theorem \ref{ClassificationApprox}.

Crudely the plan is to solve the classification problem at the von Neumann level, and lift this back to the $C^*$-setting. Slightly more precisely, we look for a quotient $\mathcal R_B$ of $B_\infty$ into which we can classify maps $A\to R_\infty$ by traces. This will fit into a short exact sequence
\begin{equation}\label{9.E1}
0\rightarrow J_B\rightarrow B_\infty\rightarrow \mathcal R_B\rightarrow 0.
\end{equation}
We then try and classify unital lifts of a given unital $^*$-homomorphism $\theta:A\to \mathcal R_B$, i.e., characterise when a lift $\tilde{\theta}:A\to B_\infty$ of $\theta$ exists, and classify these up to unitary equivalence. Successfully combining these steps will classify maps $A\to B_\infty$.  

When $B$ has unique trace $\tau$, it is natural to take $\mathcal M_B$ to be the II$_1$ factor ultrapower $(\pi_\tau(B)'')^\omega$, which is a quotient of $B_\infty$ by Kaplansky's density theorem.  Via Connes' theorem, unital $^*$-homomorphisms $\theta:A\to \mathcal R_B$ are classified up to unitary equivalence by the trace they induce on $A$ (Theorem \ref{T6.3}). Assuming additionally that $A$ has the UCT and $B$ is $\mathcal Q$--stable with $K_1(B)=0$, and working with $B_\omega$ in place of $B_\infty$ Schafhauser classified lifts of a given $\theta:A\to\mathcal R_B$ by $K_0$. Combining these two statements and then intertwining gives a classification of maps $A\to B$ by $K_0$ and traces.  Symmetrising hypotheses in the spirit of Section \ref{S6} gives the first truly abstract proof of a stably finite classification theorem. While the hypotheses are quite stringent,  they are powerful enough to show that a separable exact $C^*$-algebra satisfying the UCT and with a faithful trace embeds into a monotracial AF-algebra (\cite{S:Ann}).  It is vital that the ideal $J_B$ is separably stable, which one gets from $\mathcal Q$-stability of $B$ (the need for separable stability also forces the use of $B_\omega$ when we work with an ultrapower quotient).

Outside the unique trace setting, it is tempting to take $\mathcal R_B$ to be a suitable von Neumann ultrapower of $B^{**}_{\mathrm{fin}}$, as Connes' theorem would classify maps $A\to\mathcal R_B$ by traces.  However using positive elements $(x_n)_{n=1}^\infty$ in $B_\infty$ for which $\lim_n\tau(x_n)=0$ pointwise but not uniformly in $\tau$, one can easily obstruct separable stability of the resulting $J_B$. A more refined choice of quotient is needed to handle traces in a uniform fashion.  Such constructions came to the fore through Matui and Sato's work \cite{MS:Acta}. This is the point where Schafhauser's abstract classification machinery merges with the Toms--Winter regularity programme. Write $\|x\|_{2,T(B)}=\sup_{\tau\in T(B)}\tau(x^*x)^{1/2}$. Then we define the \emph{uniform tracial sequence algebra} by
\begin{equation}
B^\infty=\ell^\infty(B)/\{(x_n)_{n=1}^\infty\in \ell^\infty(B):\lim_{n\to\infty}\|x_n\|_{2,T(B)}=0\}.
\end{equation}
In this way, $B_\infty$ quotients onto $B^\infty$ leading to the \emph{uniform trace-kernel extension}
\begin{equation}\label{9.3}
0\to J_B\stackrel{j_b}\to B_\infty\stackrel{q_b}{\to} B^\infty\to 0.
\end{equation}
This is the right framework to classify unital maps $A\to B^\infty$ by traces and their lifts back to $B_\infty$ by the other aspects of $\underline{K}T_u$. The former uses regularity techniques, while the latter uses abstract classification. The essential point is that $\mathcal Z$-stability of $B$ gives separable stability of $J_B$ allowing $KK(A,J_B)$ to be used. In what follows, I pretend that $J_B$ is stable.

\subsection{Classifying unital maps $A\to B^\infty$}\label{S9.1}

When $B$ has a unique trace, $B^\infty$ is not quite the von Neumann algebra ultrapower $(\pi_\tau(B)'')^\omega$ used in Schafhauser's unique trace UHF-stable argument.  But for the purposes of classifying maps from separable nuclear $C^*$-algebras $A$, there is no real difference. The real challenge comes when $B$ has infinitely many extremal traces, particularly if the extreme boundary $\partial_eT(B)$ is not compact. In this case, for each $\tau\in T(B^\infty)$, Connes' theorem classifies maps $\theta:A\to \pi_\tau(B^\infty)''$. We must combine these into a classification of maps $A\to B^\infty$ by traces.

Problems of this nature have been at the heart of work on Toms--Winter regularity conjecture, and a general strategy for gluing properties from each $\pi_\tau(B^\infty)''$ together to obtain global statements which hold uniformly in all traces was developed in \cite{CETWW}.  These techniques give $B^\infty$ a ``von Neumann-like'' flavour when $B$ is nuclear and $\mathcal Z$-stable, and in particular they can be used to obtain the required classification of maps $A\to B^\infty$. Consequently, given two maps $\phi_1,\phi_2:A\to B_\infty$ which agree on traces, the compositions $q_B\circ\phi_1$ and $q_B\circ\phi_2$ are unitarily equivalent via a unitary $u\in B^\infty$ say. For each trace $\tau$ on $B^\infty$, we can write $\pi_\tau(u)$ as an exponential $e^{ih_\tau}$ for a self-adjoint $h_\tau\in \pi_\tau(B^\infty)''$. Another application of the gluing procedure can be used to find a single self-adjoint $h\in B^\infty$ with  $u=e^{ih}$. As such $u$ lifts to a unitary in $B_\infty$, and by conjugating by $u$, we can assume that $q_B\circ\phi_1=q_B\circ\phi_2$.  In this way, the remainder of the uniqueness problem for a pair of maps $\phi_1,\phi_2:A\to B_\infty$, becomes a question about the uniqueness of lifts of the common $^*$-homomorphism $q_B\circ\phi_1=q_B\circ\phi_2:A\to B^\infty$ back to $B_\infty$.

In fact, the full force of $\mathcal Z$-stability is not needed, and one can get away with a weaker central sequence condition in the spirit of Murray and von Neumann's property $\Gamma$. There is a lot going on behind the scenes here, and I will describe the ideas behind these techniques a bit further in the more regularity focused companion survey \cite{W:Survey2}.

\subsection{Classifying unital lifts}\label{S9.2}

In the second part we are given a unital map $\theta:A\to B^\infty$ and aim to classify lifts back to $B_\infty$.  A necessary condition for a lift is the existence of $\kappa\in KK(A,B_\infty)$ with $[q_B]\kappa = [\theta]$ in $KK$. We can produce these $\kappa$ using the universal (multi)coefficient theorem from the total $K$-theory component of a map $\underline{K}T_u(A)\to\underline{K}T_u(B^\infty)$.

Given such a $\kappa$, small modifications of Schafhauser's proof of the quasidiagonality theorem produces some lift $\psi_-:A\to B_\infty$ of $\theta$. But this might not have $[\psi_-]=\kappa$ in $KK(A,B_\infty)$, and it must be corrected so that it does.   By construction $\kappa-[\psi_-]$ will map to $0$ under the map $KK(A,B_\infty)\to KK(A,B^\infty)$ induced by (\ref{9.3})  so half-exactness of $KK(A,\cdot)$ gives that $\kappa-[\psi_-]$ is in the image of $KK(A,j_B):KK(A,J_B)\to KK(A,B_\infty)$.  Write $\kappa-[\psi_-]=KK(A,j_B)(\lambda)$ for some $\lambda\in KK(A,J_B)$.

Cuntz's quasihomomorphism picture of $KK$-theory is particularly well suited to $C^*$-classification problems. This defines $KK(A,J_B)$ as homotopy classes of \emph{Cuntz-pairs}: maps $(\phi_+,\phi_-):A\to\mathcal M(J_B)$ such that $\phi_+(x)-\phi_-(x)\in J_B$.  In order to translate between $KK$-theory and $^*$-homomorphisms into $B_\infty$, we need $KK$-existence and uniqueness theorems, both of which rely on absorption.  The existence theorem says that if $\phi_-:A\to \mathcal M(J_B)$ is absorbing, then given any $\lambda\in KK(A,J_B)$ we can find $\phi_+:A\to \mathcal M(J_B)$ such that the Cuntz-pair $(\phi_+,\phi_-)$ realises $\lambda$. This works in vast generality and has been regularly used in classification (in our situation all one needs is the separable-stability to work with $KK(A,J_B)$).  Following the map $\psi_-$ above by the natural map $B_\infty\to \mathcal M(J_B)$ gives rise to $\phi_-:A\to \mathcal M(J_B)$. Using the Elliott--Kucerovsky theorem (and modulo the deunitisation trick alluded to at the end of Section \ref{S8}, which is suppressed here,) $\phi_-$ is absorbing.  Thus we can find $\phi_+:A\to\mathcal M(J_B)$ such that $(\phi_+,\phi_-)$ forms a Cuntz-pair representing $\lambda$. A fairly standard pull back calculation then produces a map $\psi_+:A\to B_\infty$ also lifting $\theta$ (as a consequence of $(\phi_+,\phi_-)$ being a Cuntz-pair) so that $KK(A,j_B)(\lambda)=[\psi_+]-[\psi_-]$. Therefore $\psi_+$ realises the element $\kappa\in KK(A,B_\infty)$.  

How unique is $\psi_+$? Given two lifts $\psi_1,\psi_2:A\to B_\infty$ of $\theta$, we obtain a Cuntz-pair $(\phi_1,\phi_2):A\to \mathcal M(J_B)$ representing a class in $KK(A,J_B)$.   A $KK$- or $KL$-uniqueness theorem is designed to give asymptotic
or approximate unitary equivalence of absorbing $\phi _{1}$ and
$\phi _{2}$ when $[\phi _{1},\phi _{2}]$ vanishes in $KK(A,J_{B})$ and
$KL(A,J_{B})$, respectively. While $KK$-existence holds very generally, $KK$ and $KL$-uniqueness are more subtle, going back to Dadarlat and Eilers in the setting of $KK(A,\mathcal K)$, and it is currently unclear how generally such results can hold. For us, $\mathcal Z$-stability of $B$ is the key  ingredient through a \emph{$\mathcal Z$-stable $KL$-uniqueness theorem} developed in \cite{CGSTW} (extending a $\mathcal Q$-stable $KK$-uniqueness theorem from \cite{S:Ann}). This gives approximate unitary equivalence (with unitaries in the unitisation of $J_B\otimes\mathcal Z$) of the $\mathcal Z$-stabilisations $\phi_1\otimes 1_{\mathcal Z},\phi_2\otimes 1_{\mathcal Z}:A\to \mathcal M(J_B)\otimes\mathcal Z$ from $[\phi_1,\phi_2]=0$ in $KL(A,J_B)$.  Using separable $\mathcal Z$-stability of $B_\infty$ and the fact we work in a sequence algebra, this gives unitary equivalence of $\psi_1$ and $\psi_2$. So lifts are classified by $KL(A,J_B)$, which fits into an exact sequence
\begin{equation}
\ker KL(A,j_B)\to KL(A,J_B)\to KL(A,B_\infty).
\end{equation}
Dadarlat and Loring's universal (multi)coefficient theorem (obtained from the UCT) computes $KL(A,B_\infty)\cong \mathrm{Hom}(\underline{K}(A),\underline{K}(B_\infty))$.  We need to interpret $\ker KL(A,j_B)$ in terms of $\Ka$ and in particular the rotation maps $r_\phi$ from (\ref{7.4}).

This is achieved through an isomorphism
\begin{equation}
R_{A,B}:\ker KL(A,j_B)\to \mathrm{Hom}(K_1(A)/\mathrm{Tor}(K_1(A)),\ker\minusa_{B_\infty})
\end{equation}
with the property that $R_{A,B}([\psi_1,\psi_2])\circ t_A=r_{\psi_1}-r_{\psi_2}$ when $(\psi_1,\psi_2):A\to B_\infty$ realise a class in $\ker KL(A,j_B)$. Here $t_A:K_1(A)\to K_1(A)/\mathrm{Tor}(K_1(A))$ is the quotient map, which removes torsion from $K_1(A)$. While the individual rotation maps $r_{\psi_1}$ and $r_{\psi_2}$ depend on a choice of decomposition in (\ref{KASplit}), when $\psi_1$ and $\psi_2$ agree on $KT_u$ the difference $r_{\psi_1}-r_{\psi_2}$ does not. In this case $(r_{\psi_1}-r_{\psi_2})\circ\minusa_A=\Ka(\psi_1)-\Ka(\psi_2):\Ka(A)\to\Ka(B_\infty)$.  Then, given
$\phi _{1},\phi _{2}:A\to B_{\infty }$ agreeing on
$\underline{K}T_{u}$, by successively using traces (to reduce to the case that $q_B\circ \phi_1=q_B\circ\phi_2)$, and then total $K$-theory, and
$\overline{K}_{1}^{\mathrm{alg}}$ (to see that $\phi_1$ and $\phi_2$ induce the same class in $KL(A,J_B)$), these tools combine to give unitary equivalence of
$\phi _{1}$ and~$\phi _{2}$.

The pairing maps $\zeta^{(n)}_\bullet$ from (\ref{NewPairing}) are required for existence.  When we attempt to realise maps $\alpha:\underline{K}(A)\to\underline{K}(B_\infty)$, $\beta:\Ka(A)\to\Ka(B_\infty)$ and $\gamma:\Aff\,T(A)\to\Aff\,T(B_\infty)$ one first constructs $\theta:A\to B^\infty$ using $\gamma$. Then one lifts to $\phi:A\to B_\infty$ realising a lift $\kappa$ of $\alpha$.  As both $\Ka(\phi)$ and $(\alpha,\beta,\gamma)$ are $\underline{K}T_u$-morphisms, one can use compatibility with $\zeta^{(n)}_\bullet$ to show that the rotation map induced by $\beta-\Ka(\phi)$ vanishes on $\mathrm{Tor}(K_1(A))$.  This enables $R_{A,B}$ to be used to modify the behaviour of $\phi$ on $\Ka(A)$.

The isomorphism $R_{A,B}$ is an abstract sequence algebra version of the rotation map computations developed by Lin (see \cite{Lin:Invent}, for an example of the use of a rotation map in an asymptotic classification result). It is  built in two steps. First, the UCT gives an isomorphism $\ker(KL(A,j_B))$ to $\ker\mathrm{Hom}(\underline{K}(A),\underline{K}(j_B))$. The second part of the isomorphism is then a direct computation, which relies on the ``von Neumann like'' structure of $B^\infty$ (in particular, $K_*(B^\infty)\cong (\Aff\,T(B^\infty),0)$) via the techniques in \cite{CETWW} hinted at in Section \ref{S9.1}.

\section*{Acknowledgements}
I would like to thank my collaborators on \cite{CETWW,CGSTW} : Jos\'e Carri\'on, Jorge Castillejos, Sam Evington, Jamie Gabe, Chris Schafhauser, Aaron Tikuisis, and Wilhelm Winter. I have learnt a lot from you all.  I would like to thank Bruce Blackadar, Jos\'e Carri\'on, Jamie Gabe, Shanshan Hua, Robert Neagu, and Chris Schafhauser for their very helpful comments on an earlier draft of this article.  The paper \cite{CGSTW} was initiated at a BIRS workshop in 2017, and has benefitted massively from discussions at the American Institute of Mathematics as part of their SQuaRE programme.  I thank both institutes, and their funders for their support.


\begin{thebibliography}{99}



\bibitem{CGSTW}
J.~Carri\'on, J.~Gabe, C.~Schafhauser, A.~Tikuisis, S.~White, Classifying $^*$-homomorphisms I: simple nuclear $C^*$-algebras.  Preprint.

\bibitem{CETWW}
J.~Castillejos, S.~Evington, A.~Tikuisis, S.~White, W.~Winter, Nuclear dimension of simple $C^*$-algebras.
\emph{Invent. Math.} \textbf{224} (2021), no. 1, 245–290.

\bibitem{C:Ann}
A.~Connes, Classification of injective factors. Cases $II_1$, $II_\infty$, $III_\lambda$, $\lambda\neq 1$.
\emph{Ann. of Math. (2)} \textbf{104} (1976), no. 1, 73–115.

\bibitem{C:Invent}
E.~Christensen, Perturbations of operator algebras.
\emph{Invent. Math.} \textbf{43} (1977), no. 1, 1–13.

\bibitem{CSSWW} 
E.~Christensen, A.~Sinclair, R.~Smith, S.~White, and W.~Winter, Perturbations of nuclear $C^*$-algebras.
\emph{Acta Math.} \textbf{208} (2012), no. 1, 93–150.

\bibitem{DL}
M.~Dadarlat, T.~Loring, A universal multicoefficient theorem for the Kasparov groups.
\emph{Duke Math.} J. \textbf{84} (1996), no. 2, 355–377.

\bibitem{E:ICM} 
G.~Elliott, The classification problem for amenable $C^*$-algebras.  In \emph{Proceedings of the International Congress of Mathematicians, Vol. 1, 2 (Zürich, 1994)}, 922–932, Birkh\"auser, Basel, 1995.

\bibitem{EGLN}
G.~Elliott, G.~Gong, H.~Lin, Z.~Niu, On the classification of simple amenable $C^*$-algebras with finite decomposition rank, II
arXiv:1507.03437.

\bibitem{EK}
G.~Elliott, D.~Kucerovsky, An abstract Voiculescu--Brown--Douglas--Fillmore absorption theorem.
\emph{Pacific J. Math.} \textbf{198} (2001), no. 2, 385–409.

\bibitem{Gabe}
J.~Gabe, Classification of $\mathcal O_\infty$-stable $C^*$-algebras
\emph{ Mem. Amer. Math. Soc.}, to appear, arXiv:1910.06504.


\bibitem{GLN1} 
G.~Gong, H.~Lin, Z.~Niu, A classification of finite simple amenable $\mathcal Z$-stable $C^*$-algebras, I: $C^*$-algebras with generalized tracial rank one.
\emph{C. R. Math. Acad. Sci. Soc. R. Can.} \textbf{42} (2020), no. 3, 63–450.

\bibitem{GLN2} 
G.~Gong, H.~Lin, Z.~Niu, A classification of finite simple amenable $\mathcal Z$-stable $C^*$-algebras, II: $C^*$-algebras with rational generalized tracial rank one. 
\emph{C. R. Math. Acad. Sci. Soc. R. Can.} \textbf{42} (2020), no. 4, 451–539.

\bibitem{H:Acta}
U.~Haagerup, Connes' bicentralizer problem and uniqueness of the injective factor of type $III_1$.
\emph{Acta Math.} \textbf{158} (1987), no. 1-2, 95–148.

\bibitem{K:JEMS}
D.~Kerr, Dimension, comparison, and almost finiteness.
\emph{J. Eur. Math. Soc. (JEMS)} \textbf{22} (2020), no. 11, 3697–3745.

\bibitem{KN:arXiv}
D.~Kerr, P.~Naryshkin, Elementary amenability and almost finiteness.
arXiv:2107.05273.

\bibitem{K:ICM}
E.~Kirchberg, Exact $C^*$-algebras, tensor products, and the classification of purely infinite algebras. In \emph{Proceedings of the International Congress of Mathematicians, Vol. 1, 2 (Zürich, 1994)}, 943–954, Birkhäuser, Basel, 1995.

\bibitem{KP:Crelle}
E.~Kirchberg, N.C.~Phillips, Embedding of exact $C^*$-algebras in the Cuntz algebra $\mathcal O_2$.
\emph{J. Reine Angew. Math.} 525 (2000), 17–53.

\bibitem{L:Invent} 
X.~Li, Every classifiable simple $C^*$-algebra has a Cartan subalgebra.
\emph{Invent. Math.} \textbf{219} (2020), no. 2, 653–699.

\bibitem{Lin:Invent}
H.~Lin, Asymptotic unitary equivalence and classification of simple amenable $C^*$-algebras.
\emph{Invent. Math.} \textbf{183} (2011), no.~2, 385--450.

\bibitem{MS:Acta}
H.~Matui, Y.~Sato, Strict comparison and $\mathcal Z$-absorption of nuclear $C^*$-algebras.
\emph{Acta Math.} \textbf{209} (2012), no. 1, 179–196.

\bibitem{MS:DMJ} 
H.~Matui, Y.~Sato, Decomposition rank of UHF-absorbing $C^*$-algebras.
\emph{Duke Math. J.} \textbf{163} (2014), no. 14, 2687–2708.

\bibitem{McDuff} 
D.~McDuff, Central sequences and the hyperfinite factor.
\emph{Proc. London Math. Soc. (3)} \textbf{21} (1970), 443–461.

\bibitem{MvN.4}
F.~Murray, and J.~von Neumann, On rings of operators. IV. 
\emph{Ann. of Math. (2)} \textbf{44} (1943), 716–808.


\bibitem{P:Doc}
N.C.~Phillips, A classification theorem for nuclear purely infinite simple $C^*$-algebras. 
\emph{Doc. Math.} \textbf{5} (2000), 49–114.

\bibitem{R:Acta}
M.~Rørdam, 
A simple $C^*$--algebra with a finite and an infinite projection.
\emph{Acta Math.} \textbf{191} (2003), no. 1, 109–142.

\bibitem{RS}
J.~Rosenberg, C.~Schochet, The K\"unneth theorem and the universal coefficient theorem for Kasparov's generalized $K$-functor.
\emph{Duke Math. J.} \textbf{55} (1987), no. 2, 431–474.

\bibitem{S:Crelle}
C.~Schafhauser, A new proof of the Tikuisis--White--Winter theorem.
\emph{J. Reine Angew. Math.} \textbf{759} (2020), 291–304.

\bibitem{S:Ann}
C.~Schafhauser, Subalgebras of simple AF-algebras.
\emph{Ann. of Math. (2)} \textbf{192} (2020), no. 2, 309–352.

\bibitem{Thomsen}
K.~Thomsen, Traces, unitary characters and crossed products by $\mathbb Z$.
\emph{Publ. Res. Inst. Math. Sci.} \textbf{31} (1995), no. 6, 1011–1029.

\bibitem{TWW}
A.~Tikuisis, S.~White, W.~Winter, Quasidiagonality of nuclear $C^*$-algebras.
\emph{Ann. of Math. (2)} \textbf{185} (2017), no. 1, 229–284.

\bibitem{T:Ann}
A.~Toms, On the classification problem for nuclear $C^*$-algebras. 
\emph{Ann. of Math. (2)} \textbf{167} (2008), no. 3, 1029–1044.

\bibitem{W:Survey2}
S.~White, $\mathcal Z$-stability, tracial completions and regularity for $C^*$-algebras.
Manuscript in preparation.

\bibitem{W:Invent} 
W.~Winter, Nuclear dimension and $\mathcal Z$-stability of pure $C^*$-algebras. 
\emph{Invent. Math.} \textbf{187} (2012), no. 2, 259–342.

\bibitem{W:Crelle}
W.~Winter, Localizing the Elliott conjecture at strongly self-absorbing $C^*$-algebras.
J. Reine Angew. Math. \textbf{692} (2014), 193–231.

\bibitem{W:ICM}
W.~Winter, Structure of nuclear $C*$-algebras: from quasidiagonality to classification and back again. In \emph{Proceedings of the International Congress of Mathematicians—Rio de Janeiro 2018. Vol. III. Invited lectures}, 1801–1823, World Sci. Publ., Hackensack, NJ, 2018.








\end{thebibliography}
\end{document}